\documentclass[aos,MSNbibl,nameyear,dvips]{arximspdf}

\doi{10.1214/14-AOS1247} 
\volume{42}
\issue{6}
\pubyear{2014}
\firstpage{2282}
\lastpage{2300}
\docsubty{FLA}

\makeatletter
\newcommand{\lleft}{\left}
\newcommand{\rright}{\right}
\newtheorem{theorem}{Theorem}
\newtheorem{proposition}[theorem]{Proposition}
\newtheorem{lemma}[theorem]{Lemma}

\newcommand{\lset}{\mathcal{L}_K} 
\newcommand{\eset}{\mathcal{E}} 
\newcommand{\lsett}{\mathcal{L}_{\widetilde K}} 
\newcommand{\ecset}{\mathcal{C}} 
\makeatother

\begin{document}
\begin{frontmatter}

\title{Optimal cross-over designs for full interaction~models}
\runtitle{Cross-over designs}

\begin{aug}
\author[A]{\fnms{R.~A.}~\snm{Bailey}\ead[label=e1]{rab24@st-andrews.ac.uk}}
\and
\author[B]{\fnms{P.}~\snm{Druilhet}\corref{}\ead[label=e2]{pierre.druilhet@univ-bpclermont.fr}}
\affiliation{Queen Mary University of London and University of St
Andrews, and Clermont Universit\'e}
\address[A]{School of Mathematical Sciences\\
Queen Mary University of London\\
Mile End Road, London E1 4NS\\
United Kingdom\\
and\\
School of Mathematics and Statistics\\
University of St Andrews\\
St Andrews, Fife, KY16 9SS\\
United Kingdom\\
\printead{e1}}
\runauthor{R. A. Bailey and P. Druilhet}
\address[B]{Laboratoire de Math\'{e}matiques\\
UMR CNRS 6620\\
Clermont Universit\'e\\
and\\
Universit\'{e} Blaise Pascal\\
63177 AUBIERE Cedex\\
France\\
\printead{e2}}
\end{aug}

\received{\smonth{6} \syear{2013}}
\revised{\smonth{5} \syear{2014}}

%
\begin{abstract}
We consider repeated measurement designs when a residual or carry-over effect
may be present in at most one later period. Since assuming an additive
model may be unrealistic for some applications
and leads to biased estimation of
treatment effects, we consider a model with interactions between
carry-over and direct treatment
effects. When the aim of the experiment is to study the effects of a
treatment used alone, we obtain universally optimal approximate
designs. We
also propose some efficient designs with a reduced number of subjects.
\end{abstract}

%

\end{frontmatter}

\section{Introduction}\label{sec1}

In repeated measurement designs or crossover designs, interference
is often observed between a direct treatment effect and the treatment
applied in the previous period.
We denote by $\xi_{uv}$ 
the effect of treatment $u$ 
when it is preceded by treatment $v$. 
There are several ways to model such effects. The simplest one is to
assume that there is no inter\-ference. In this case, $\xi_{uv} = \tau_u$,
the direct treatment effect.

For a parsimonious interference model, we may assume that the direct
and the
carry-over effects are additive. In this case, 
$\xi_{uv}=\tau_u+\lambda_v$, where $\tau_u$ 
is the direct effect of treatment $u$, 
and $\lambda_v$ 
is the carry-over effect due to treatment $v$. 
In practice, this model is often unrealistic.

Kempton, Ferris and
David (\citeyear{kemp:ferre:davi:2001}) propose an interference model in which a treatment
which has a large direct effect will also have a large carry-over effect.
More precisely, they assume that the carry-over effect is proportional to
the direct effect. \citet{bail:kune:2006} obtain optimal designs under
this model.

\citet{afsa:heda:2002} propose another way to enrich the additive
models: they assume that the carry-over
effect of a treatment depends on whether that treatment is preceded by
itself or not. In that case $\xi_{uv}=\tau_u+\lambda_v+\chi_{uv}$,
where $\chi_{uv}=0$ if $u\neq v$ 
and $\chi_{uu}$ 
represents the specific effect of treatment $u$ 
preceded by itself. For this model, optimal designs are obtained by
Kunert and Stufken (\citeyear{kune:stuf:2002,kune:stuf:2008}) when the parameters of interest are the
direct treatment effects, and by \citet{drui:tins:2009} when the
parameters of interest are the total effects $\tau_u + \lambda_u + \chi_{uu}$.

The finest possible model, proposed by \citet{sen:muke:1987}, assumes
full interactions between carry-over
and direct treatment effects, which means that no constraints on $\xi_{uv}$
are assumed. For a full interaction model, there is no natural way to
define a direct treatment effect. For example, \citet{park:bose:all:2011}
obtain efficient designs when the parameters of interest are the
standard least-squares means of treatments, that is, $t^{-1}\sum_v\xi
_{uv}$ for $1\leq u\leq t$,
where $t$ is the number of treatments to be compared. Under a full
inter\-action model, the contrasts of the least-squares means depend on
all the other treatment effects through their interactions.

When the aim of the experiment is to select a single treatment which
will be used alone, that is, preceded by itself, the relevant effects
to be considered are total effects $\phi_u = \xi_{uu}$ for $1\leq u
\leq t$,
which correspond to the effect of a treatment preceded by itself; see
\citet{bail:drui:2004} for a review of situations where total
effects have to be considered.

\citet{kush:1997} and \citet{kune:mart:2000} propose a method for
obtaining optimal cross-over designs for direct treatment effects in
the framework of approximate designs by using Schur-complement properties.
The method has three main steps: (i) expressing the information matrix
of the whole design as a sum of the information matrices for the
sequences of treatments given to individual subjects (Section~\ref
{subsection.Cphi});
(ii) considering so-called symmetric designs, in which the proportion
of subjects given any sequence is invariant under the symmetric group
of all permutations of the treatments (Section~\ref{subsection.symmetric});
applying maximin procedures to equivalence classes of sequences
(Section~\ref{section.universal}).

A first generalisation of these techniques for more general effects is
proposed by \citet{drui:tins:2009}.
In this paper, we propose a higher level of generalisation by using
group theory to obtain optimal designs for total effects under the full
interaction interference model. We also propose efficient designs of
reduced sizes.

\section{The designs and the model}
\label{section.designandmodel}

We consider a design $d$ with $n$ subjects and~$k$ periods. Let $t$ be the
number of treatments. For $1\leq i\leq n$ and $1\leq j\leq k$, denote by
$d ( i,j ) $ the treatment assigned to subject $i$ in period
$j$. We
assume the following full treatment $\times$ carry-over interaction
model for the response $y_{ij}$:
%
\begin{equation}
\label{eq.model.sansperiode} y_{ij}=\beta_{i}+ \xi_{d(i,j),d(i,j-1)}+
\varepsilon_{ij},
\end{equation}
where $\beta_{i}$ is the effect of subject $i$, and $\xi_{uv}$ is the
effect of treatment $u$ when preceded by treatment $v$. For the first
period, we assume a specific carry-over effect that can be represented
by a fictitious treatment labelled $0$: $\xi_{u0}$ represents the
effect of treatment $u$ with no treatment before. The residual errors
$\varepsilon_{ij}$ are assumed to be independent and identically
distributed with
expectation $0$ and variance $\sigma^{2}$.
In most applications, a period effect is included in the model. It will
be seen in Section~\ref{subsection.periode}
that optimal designs found for model (\ref{eq.model.sansperiode}) are
also optimal when period effects are added.

In vector notation, model (\ref{eq.model.sansperiode}) can be written
\[
Y=B\beta+ X_{d}\xi+\varepsilon,
\]
where $Y$ is the $nk$-vector of responses with entries $y_{ij}$ in
lexicographic order,
and $\beta$ is the $n$-vector of subject effects.
The entries of the $t(t+1)$-vector $\xi$ are denoted by $\xi_{uv}$ and
sorted in lexicographic order.
The matrices associated with these effects are, respectively, given by $B
$ and $X_{d}$.
Note that $B=I_n\otimes\mathbb I_k$, where $I_n$ denotes the identity
matrix of order $n$, the symbol $\otimes$
denotes the Kronecker product, and $\mathbb I_k$ is the
$k$-dimensional vector of ones.
Also, $X_d$ is an $nk \times t(t+1)$ matrix whose entries are all $0$
apart from a single $1$ in each row. In particular, $X_d \mathbb
{I}_{t(t+1)} = \mathbb{I}_{nk}$.
We have $\mathbb E(\varepsilon)=0$ and $\operatorname{Var}%
 ( \varepsilon ) =\sigma^{2}I_{nk}$.

We denote by $\phi$ the $t$-vector of total effects, which corresponds
to the situation where a treatment is preceded by itself. We have
$\phi_u=\xi_{uu}$, for $u=1, \ldots, t$. 
Denote by $ K$ the $t(t+1)\times t$ matrix with entries $K_{uv}^w=1$ if
$u=v=w$ and $0$ otherwise for $u,w=1,\ldots,t$ and $v=0,\ldots,t$,
where $w$ is the single index for the columns, and $uv$ is
the double index for the rows, similar to the index for the vector
$\xi_{uv}$. We have
%
\begin{equation}
\label{eq.definitionphi} \phi=K^{\prime}\xi.
\end{equation}

\section{Information matrices for total effects}

\subsection{Information matrix for \texorpdfstring{$\xi$}{xi} and \texorpdfstring{$\phi$}{phi}}
\label{sub:infomat}

Put $\omega_{B}=B ( B^{\prime}B ) ^{-1}B'$, which is the projection
matrix onto the column space\vspace*{1pt} of $B$, and
$\omega_{B}^{\bot}=I_{nk}-\omega_{B}=I_{n}\otimes Q_{k}$ with
$Q_{k}=\omega_{\mathbb{I}_{k}}^{\bot}=I_{k}-k^{-1} J_k$,
where $ J_k=\mathbb I_{k}\mathbb I_{k}'$.
The information matrix $C_d [ \xi ] $ for the vector $\xi$ is
given by [see, e.g., \citet{kune:1983}]
\[
C_{d} [ \xi ] =X_d' \omega_{B}^{\bot}
X_d.
\]
Note that $\omega_B^\bot X_d \mathbb{I}_{t(t+1)} = \omega_B^\bot\mathbb
{I}_{nk} = \mathbf{0}$, and so
%
\begin{equation}
C_{d} [ \xi ] \mathbb{I}_{t(t+1)} =\mathbf{0}. \label{eq:zerosum}
\end{equation}

Denote by $X_{di}$ the $k\times t(t+1)$ design matrix for subject $i$ and
by $C_{di}[\xi]=X_{di}' Q_k X_{di}$ the information matrix
corresponding to subject $i$ alone. We have $X_d'=(X_{d1}',\ldots,X_{dn}')$ and
\[
C_{d} [ \xi ] =\sum_{i=1}^n
C_{di}[\xi]= \sum_{i=1}^n
X_{di}' Q_k X_{di}. %
\]

Note that $X_{di}$ and therefore $C_{di}[\xi]$ depend only on the
sequence of treatments applied to subject $i$.
Denote by $\mathcal S$ the set of all sequences of $k$ treatments.
For a design $d$ and a sequence $s\in\mathcal S$, denote by $\pi_d(s)$
the proportion of subjects that receive $s$,
and denote by $X_s$ and $C_{s}[\xi]$ the associated matrices. We have
%
\begin{equation}
\label{eq.infomatxi} C_{d} [ \xi ] =n \sum_{s\in\mathcal S}
\pi_d(s) C_{s}[\xi]= n \sum_{s\in\mathcal S}
\pi_d(s) X_{s}' Q_k
X_{s}.
\end{equation}

\label{subsection.Cphi}
The information matrix for the parameter of interest
$\phi=K'\xi$ may be obtained from $C_d[\xi]$ by the extremal
representation [see
\citet{gaff:1987} or \citet{puke:1993}]
%
\begin{equation}
\label{eq.infomatphi} C_{d} [ \phi ] =C_{d} \bigl[ K^{\prime}
\xi \bigr] =
\min_{L \in\lset}
L^{\prime}C_{d} [\xi ] L,
\end{equation}
where $\lset= \{L\in\mathbb{R}^{t(t+1)\times t} \mid L^{\prime
}K=I_{t}\}$
and the minimum is taken relative to the Loewner ordering. The minimum in
(\ref{eq.infomatphi}) exists and is unique for a given design $d$.
Put $\eset_d = \{L\in\lset\mid L'C_d[\xi]L = C_d[\phi]\}$.

In the sequel, the entries of $L$, or, more generally, of any matrix
of size $t(t+1)\times t$, will be denoted by $L^w_{uv}$, for
$u,w=1,\ldots,t$, and $v=0,\ldots,t$, where $w$ is the column index and
$uv$ is the double index for the rows, similar to the vector $\xi$ or
the matrix $K$. The $t\times t$ matrix $L'K$ has entries
$(L'K)_{uv}=L^u_{vv}$, for $u,v=1,\ldots,t$.

\begin{lemma}
\label{lemma.Cdphi1=0}
For any design $d$, the row and column sums of $C_d[\phi]$ are zero.
\end{lemma}

\begin{pf}
Since $C_d[\phi]$ is symmetric, we have to prove that $\mathbb
I_t'C_d[\phi]\mathbb I_t=0$.
Consider the $t(t+1)\times t$ matrix $L$ such that $L^u_{vw}$ is equal
to $1$ if $u=v$ and $0$ otherwise.
The matrix $L$ satisfies $L \mathbb I_t=\mathbb I_{t(t+1)}$ and the
constraint $L'K=I_t$.
It follows from~(\ref{eq.infomatphi})
and (\ref{eq:zerosum})
that
$ 0 \leq\mathbb I_t' C_d[\phi] \mathbb I_t\leq\mathbb I_t' L^{\prime}
C_{d} [
\xi ] L \mathbb I_t= \mathbb I_{t(t+1)}' C_{d} [\xi ]
\mathbb I_{t(t+1)} = 0$.
\end{pf}

For a design $d$, denote by $L^{\ast}$
a matrix in $\eset_d$.
Since, for any given $L$, $L^{\prime}C_{d} [ \xi ] L$ is
linear in $C_{d} [\xi ]$,
we have by (\ref{eq.infomatxi}),
%
\begin{equation}
\label{eq.cphiLstar} C_d[\phi] =L^{\ast\prime} C_{d} [ \xi ]
L^{\ast}= n\sum_{s\in\mathcal S} \pi_d(s)
L^{\ast\prime} C_{s}[\xi] L^{\ast}.
\end{equation}
This linearisation is the basis of Kushner's methods.

\subsection{Approximate designs and symmetric designs}
\label{subsection.symmetric}

An exact design is characterised, up to a subject permutation,\vspace*{1pt} by the
proportions of sequences that appear in it. 
These proportions are multiples of $n^{-1}$. 
If we allow the proportions to vary continuously in $[0,1]$ with the
only restriction that the sum must be equal to $1$,
we obtain an approximate design. By definition, the information
matrices of $\xi$ and $\phi$ for an approximate designs are given by
(\ref{eq.infomatxi})
and (\ref{eq.infomatphi}) as for an exact design.
The second idea of Kushner's method is to find a universally optimal
design in the set of approximate designs using the linearised
expression (\ref{eq.cphiLstar}).
If the optimal approximate design is not an exact design, one can
calculate a sharp lower bound for efficiency factors of competing exact designs.

We now recall the concepts of permuted sequence, symmetric design, and
symmetrised design as introduced by \citet{kush:1997}.
Let $\sigma$ be a permutation of the treatment labels $\{1,\ldots,t\}$
and $s$ a sequence of treatments. The \emph{permuted sequence} $s_\sigma
$ is obtained from $s$ by permuting the treatment labels according to
$\sigma$. Similarly, the design $d_\sigma$ is the design obtained from
the design $d$ by permuting the treatment labels according to $\sigma$.
A design $d$ is said to be a \emph{symmetric design} if, for any
sequence $s$ and any permutation $\sigma$, $\pi_d(s_\sigma)=\pi_d(s)$.
For such a design, $d$ and $d_\sigma$ are identical up to a subject
permutation, which may be written $d=d_\sigma$. From a design~$d$, we
define the \emph{symmetrised design} $\bar d$ by
%
\begin{equation}
\label{eq.pisymmetrised} \pi_{\bar d}(s)=\frac{1}{t!}\sum
_{\sigma\in S_t}\pi_d(s_{\sigma})\qquad \forall s\in
\mathcal S,
\end{equation}
where $S_t$ is the set of all permutations of $\{1,\ldots,t\}$. It is
easy to see that the symmetrised design $\bar d$ is a symmetric design.

To a permutation $\sigma$ of treatment labels, we may associate a
permutation $\sigma^*$ of the carry-over effect labels $\{0,1,\ldots,t\}$
where $\sigma^*(0)=0$ and $\sigma^*(u)=\sigma(u)$ for $u=1,\ldots,t$.
We also associate a permutation $\widetilde\sigma$ of
$\{1,\ldots,t\}\times\{0,\ldots,t\}$ defined by
$\widetilde\sigma(u,v)=(\sigma(u),\sigma^*(v))$. We denote by $P_\sigma$,
$P_{\sigma^*}$, and $P_{\widetilde\sigma}=P_\sigma\otimes P_{\sigma^*}$ the
corresponding permutation matrices:
for example, $P_\sigma(u,v)=1$ if $\sigma(u)=v$ and $P_\sigma(u,v)=0$ otherwise.

For $L\in\lset$, put $ L_\sigma= P_{\widetilde\sigma}' LP_\sigma$.
It can be checked that $ P_{\widetilde\sigma} ' KP_\sigma=K$; see also
the definition of the matrix $L_{(1)}$ after Lemma~\ref{lemma.Lstarsymmetric}.

\begin{lemma}
\label{lem:alldes}
For any design $d$ and any permutation $\sigma$ in $S_t$, we have:
%
\begin{eqnarray}
\label{eq.Cxidsigma} C_{d_{\sigma}}[\xi] &=& P_{\widetilde\sigma}
C_d[\xi] P_{\widetilde\sigma}';
\\
\label{eq.Cphidsigma} C_{ d_{\sigma}}[\phi] &=& P_{\sigma} C_d[
\phi] P_{\sigma}';
\\
\label{eq.Cxidbar} C_{\bar d}[\xi] &=& \frac{1}{t!} \sum
_{\sigma\in
S_t} P_{\widetilde\sigma} C_d[\xi]
P_{\widetilde\sigma}';
\\
\label{eq.Cphidbar} C_{\bar d}[\phi] &\geq& \frac{1}{t!} \sum
_{\sigma
\in S_t} P_{\sigma} C_d[\phi]
P_{\sigma}' \qquad\mbox{w.r.t. the Loewner ordering};
\end{eqnarray}
and $L\in\eset_{d}$ if and only if $L_\sigma\in\eset_{d_\sigma}$.
\end{lemma}

\begin{pf}
By definition of $P_{\widetilde\sigma}$, $X_{d_\sigma}=X_dP_{\widetilde
\sigma}' $, and so $C_{d_{\sigma}}[\xi]=X'_{d_{\sigma}} \omega_{B}^{\bot
} X_{d_{\sigma}} =
P_{\widetilde\sigma}X'_{d} \omega_{B}^{\bot} X_{d} P_{\widetilde\sigma
}'=P_{\widetilde\sigma} C_d[\xi] P_{\widetilde\sigma}' $, which
corresponds to
(\ref{eq.Cxidsigma}).
If $L\in\lset$, then\break $L'C_{d_\sigma}[\xi]L =
L' P_{\widetilde\sigma} C_d[\xi] P_{\widetilde\sigma}' L =
P_\sigma L_\sigma' C_d[\xi] L_\sigma P_\sigma'$.
Now $L_\sigma'K=
P_\sigma' L'P_{\widetilde\sigma}P_{\widetilde\sigma}' KP_\sigma=
P_\sigma' L 'KP_\sigma$.
If $L\in\lset$, then $L'K =I_t$, so $L_\sigma'K=I_t$ and $L_\sigma\in
\lset$.
The same argument with $\sigma^{-1}$ shows that if $L_\sigma\in\lset$
then $L\in\lset$.
The Loewner ordering is unchanged by permutations, so
\[
C_{ d_{\sigma}}[\phi] = \min_{L\in\lset} \bigl(L'
C_{d_\sigma}[\xi ]L \bigr) = P_\sigma \Bigl(\min
_{L_\sigma\in\lset} L_\sigma'C_d[
\xi]L_\sigma \Bigr) P_\sigma' =P_\sigma
C_{d}[\phi] P_\sigma',
\]
and (\ref{eq.Cphidsigma}) is established.
Moreover, $L\in\eset_d$ if and only if $L_\sigma\in\eset_{d_\sigma}$.
Formula (\ref{eq.Cxidbar}) follows directly from (\ref
{eq.Cxidsigma}) and (\ref{eq.pisymmetrised}). Formula (\ref
{eq.Cphidbar}) follows from (\ref{eq.Cxidbar}) and the concavity of
the minimum representation
(\ref{eq.infomatphi}).
\end{pf}

We recall that a $t\times t$ matrix $C$ is completely symmetric if $C=a
I_t + b J_t$ for some scalars $a$ and $b$ or, equivalently, if
$P_{\sigma} C P_{\sigma}'=C$ for every permutation $\sigma$ in~$S_t$.

\begin{lemma}
\label{lemma.cphics}
If $d$ is a symmetric design, then $C_d[\phi]$ is completely symmetric.
\end{lemma}

\begin{pf}
Since $d$ is symmetric, $d_\sigma=d$. By (\ref{eq.Cphidsigma}),
$C_d[\phi]=C_{d_\sigma}[\phi]=P_{\sigma} C_d[\phi] P_{\sigma}'$ for any
permutation $\sigma$ in $S_t$. Therefore $C_d[\phi]$ is completely symmetric.
\end{pf}

The key point to obtain an optimal design is to identify the structure
of the $t(t+1)\times t$ matrix $L^*$ defined in (\ref{eq.cphiLstar}),
whose entries are denoted by $L^{*w}_{uv}$.

\begin{lemma}
\label{lemma.Lstarsymmetric}
If $d$ is a symmetric design, then the matrix $L^*$ 
in (\ref{eq.cphiLstar}) can be chosen so that it satisfies
%
\begin{equation}
\label{eq.invarianLstar} 
L_\sigma^* =L^*\qquad\forall\sigma\in
S_t,
\end{equation}
or, equivalently,
%
\begin{equation}
\label{eq.invariantorbite} L^{*\sigma(w)}_{\sigma(u)\sigma^*(v)}=L^{*w}_{uv}\qquad
\forall\sigma\in S_t.
\end{equation}
\end{lemma}

\begin{pf}
If $\sigma\in S_t$, then $d_\sigma= d$, so $\eset_{d_\sigma} = \eset
_d$, and
Lemma~\ref{lem:alldes} shows that $L_\sigma\in\eset_d$.
Put $L^*=  (\sum_{\sigma\in S_t} L_\sigma )/t!$,
which satisfies (\ref{eq.invarianLstar}). Since $\eset_d$ is closed
under taking averages [see \citet{drui:tins:2009}, proof of Lemma
A1], $L^*$ also belongs to $\eset_d$.
\end{pf}

A consequence of (\ref{eq.invariantorbite}) is that the entries
$L^{*w}_{uv}$ are constant for $(u,v,w)$ belonging to the same orbit of
the permutation group $\{(\widetilde\sigma,\sigma)\}_{\sigma\in S_t}$
acting on $\{1,\ldots,t\}\times\{0,\ldots,t\}\times\{1,\ldots,t\}$.
There are seven distinct orbits:
\begin{itemize}
%
\item$\mathcal O_1=\{(u,u,u) \mid u=1,\ldots,t\}$,
\item$\mathcal O_2=\{(u,v,u) \mid u,v=1,\ldots,t, u\neq v\}$,
\item$\mathcal O_3=\{(u,v,v) \mid u,v=1,\ldots,t, u\neq v\}$,
\item$\mathcal O_4=\{(u,v,w) \mid u,v,w=1,\ldots,t, u\neq v\neq w
\neq u\}$,
\item$\mathcal O_5=\{(u,0,u) \mid u=1,\ldots,t\}$,
\item$\mathcal O_6=\{(u,0,w) \mid u,w=1,\ldots,t, u\neq w\}$,
\item$\mathcal O_7=\{(u,u,w) \mid u,w=1,\ldots,t, u\neq w\}$.
\end{itemize}
For $q=1,\ldots,7$, denote by $L_{(q)}$ the $t(t+1)\times t $ matrix
with entries $L_{(q)uv}^w=1$ if $(u,v,w)$ belongs to the orbit
$\mathcal O_q$ and $0$ otherwise. Note that $L_{(1)}=K$.

By construction of $L_{(q)}$, we have
%
\begin{equation}
\label{eq.invarianceLi} P_{\widetilde\sigma}' L_{(q)}
P_\sigma=L_{(q)} \qquad\forall\sigma\in S_t\mbox{ and
}q=1,\ldots,7.
\end{equation}

\begin{proposition}
For a symmetric design $d$, the matrix $L^*$
in Lemma~\ref{lemma.Lstarsymmetric}
may be written as 
%
\begin{equation}
\label{eq.Lstarreduc} L^* = L_\gamma=L_{(1)}+\sum
_{q=2}^6\gamma_q L_{(q)},
\end{equation}
where $\gamma=(\gamma_2,\ldots,\gamma_7)$ is a vector of scalars.
\end{proposition}

\begin{pf}
Since $L^*$ satisfies (\ref{eq.invarianLstar}), it
is a linear combination of the matrices $L_{(q)}$:
$
L^*=\sum_{q=1}^7\gamma_q L_{(q)}$.
It can be checked that $L_{(1)}'K=K'K=I_t$, $L_{(7)}'K=J_t-I_t$
and $L_{(q)}'K=0$ for $q=2,\ldots,6$.
Consequently, the constraint $L^{*\prime}K=I_t$ may be written
$\gamma_1=1$ and $\gamma_7=0$.
\end{pf}

\subsection{The model with period effects}
\label{subsection.periode}
We consider here the same model as in Section~\ref
{section.designandmodel} with the addition of a period effect.
The response for subject $i$ in period $j$ is given by
%
\begin{equation}
\label{eq.modele.avecperiode} y_{ij}=\alpha_j+\beta_{i}+
\xi_{d(i,j),d(i,j-1)}+\varepsilon_{ij},
\end{equation}
where $\alpha_j$ is the effect of period $j$. In vector notation, we have
\[
Y= A\alpha+B\beta+ X_{d}\xi+\varepsilon,
\]
with $A=\mathbb I_n\otimes I_k$,
where $\alpha$ is the $k$-vector of period effects.
Denote $\theta'=(\xi',\alpha')$. If $d$ is an exact design, the
information matrix for $\theta$ is given by
\[
\widetilde C_d[ \theta]= \pmatrix{
C_d[\xi] & C_{d12}
\vspace*{2pt}\cr
C_{d21} & C_{d22}} = \pmatrix{
 X'_d \omega_{B}^{\bot}
X_d & X_d' \omega_{B}^\perp
A
\vspace*{2pt}\cr
A' \omega_{B}^\perp X_d &
A' \omega_{B}^\perp A}, %
\]
where $C_d[\xi]$ is the information matrix for $\xi$ obtained in the
model without period effects and $C_{d22}= n Q_k$.

The $t$-vector $\phi$ of total effects defined by (\ref
{eq.definitionphi}) may also be seen as a subsystem of the parameter
$\theta$, 
because $\phi=\widetilde K'\theta$ with $\widetilde K'=(K',0_{t\times k})$.
The information matrix $\widetilde C_d[\phi]$ for $\phi$ under model
(\ref{eq.modele.avecperiode}) may be obtained from $\widetilde
C_d[\theta]$ by the extremal representation
\[
\label{eq.infomatphitilde} \widetilde C_{d} [ \phi ] = 
\min
_{\widetilde{L} \in\lsett} 
\widetilde
L^{\prime}C_{d} [ \theta ] \widetilde L, %
\]
where $\lsett= \{\widetilde{L}\in\mathbb{R}^{(t(t+1)+k)\times t} \mid
\widetilde L^{\prime}\widetilde K=I_{t}\}$. Partitioning $\widetilde
L'$ as $(L'\mid N')$ with $L$ and $N$ of
sizes $t(t+1)\times t$ and $k\times t$, we have
%
\begin{equation}
\label{eq.CminLperiode} 
\widetilde C_{d} [ \phi ] = \min
_{(L' \mid N')' \in\lsett} 
 \bigl( L^{\prime}C_{d}
[ \xi ] L +L'C_{d12}N+N'C_{d21}L+N'C_{d22}N
\bigr).
\end{equation}
Note that $(L' \mid N')' \in\lsett$
is equivalent to $L\in\lset$ for $L$ and $N$ with suitable dimensions.
Choosing $N=0$ in (\ref{eq.CminLperiode}), we have
$
\widetilde C_d[\phi]\leq C_d[\phi]
$
with respect to the Loewner ordering,
where $C_d[\phi]$ is the information matrix for $\phi$ under the model
without period effects, as defined in (\ref{eq.infomatphi}).
Therefore $0\leq\mathbb I_t' \widetilde C_d[\phi] \mathbb I_t \leq
\mathbb I_t' C_d[\phi] \mathbb I_t= 0$.
Hence the row and column sums of $\widetilde C_d[\phi]$ are all zero,
and so
$Q_t \widetilde C_d[\phi] Q_t = \widetilde C_d[\phi]$.

For $\sigma\in S_t$, define the permutation $\bar\sigma$ for the entries
of $\theta$ such that the entries of $\xi$ are permuted according to
$\widetilde\sigma$ and those of $\alpha$ remain unchanged. The associated
permutation matrix $P_{\bar\sigma}$ is the block diagonal matrix with
diagonal blocks $P_{\widetilde\sigma}$ and~$I_k$.
For $\widetilde L$ in $\lsett$, put
$\widetilde L_\sigma= P'_{\bar\sigma} \widetilde L P_\sigma$.
If $\widetilde L' = (L' \mid N')$, then
$\widetilde L'_\sigma= (L'_\sigma\mid N'_\sigma)$,
where $N_\sigma= NP_\sigma$.

\begin{lemma}
For any design $d$ and any permutation $\sigma$ of treatment labels, we have
%
\begin{eqnarray}
\label{eq.C12dsigma} C_{d_{\sigma}12} &=& P_{\widetilde\sigma}
C_{d12};
\\
\label{eq.Ctildephidsigma} \widetilde C_{ d_{\sigma}}[\phi] &=& P_{\sigma} \widetilde
C_d[\phi] P_{\sigma}'.
\end{eqnarray}
\end{lemma}

\begin{pf}
Equation (\ref{eq.C12dsigma}) follows from the fact that
$X_{d_\sigma}=X_d P_{\widetilde\sigma}'$.
The proof of~(\ref{eq.Ctildephidsigma}) is similar to the proof of
(\ref{eq.Cphidsigma}),
replacing $\xi$, $L$, $\lset$, and $K$ by $\theta$, $\widetilde L$,
$\lsett$, and $\widetilde K$, respectively.
\end{pf}

%

An exact design is said to be \emph{strongly balanced on the periods}
if it satisfies the following conditions:
\begin{longlist}[(iii)]
\item[(i)] for the first period, each treatment appears equally often;
\item[(ii)] for any given period, except the first one, each treatment
appears preceded by itself equally often;
\item[(iii)] for any given period, except the first one, the number of times a
treatment, say $u$, is preceded by another treatment $v$ does not
depend on $u$ or $v$.
\end{longlist}

Note that a symmetric exact design is strongly balanced on the periods.

\begin{lemma}
\label{lemma:strong}
If a design $d$ is strongly balanced on the periods and $\sigma\in S_t$,
then $P_{\widetilde{\sigma}}' X_d'A = X_d'A$.
\end{lemma}

\begin{pf}
The $(uv,j)$-entry of $X_d'A$ is equal to the number of times that
treatment $u$ occurs in period $j$ preceded by treatment $v$. Strong
balance implies that there is a single value for $v=0$, another single
value for $v=u$, and another single value for $v\notin\{0,u\}$.
Permutation of the treatments does not change this.
\end{pf}

Given a design $d$, let $G_d$ be the subgroup of $S_t$ consisting of
those permutations $\sigma$ satisfying $d_\sigma= d$ (up to a subject
permutation).
Note that a symmetric design may be characterised by $G_d=S_t$. The
subgroup $G_d$ is said to be \emph{transitive} on $\{1, \ldots, t\}$, if,
given $u$, $v$ in $\{1, \ldots, t\}$, there is some $\sigma$ in $G_d$ with
$\sigma(u)=v$. The subgroup $G_d$ is \emph{doubly transitive} if, given
$u_1$, $u_2$, $v_1$, $v_2$ with $u_1\ne u_2$ and $v_1 \ne v_2$ there is
some $\sigma$ in $G_d$ with $\sigma(u_1)=v_1$ and $\sigma(u_2) = v_2$.

\begin{proposition}
\label{prop.equalinfomatperiode}
If $d$ is an exact design with strong balance on the periods and with
transitive group $G_d$, then the information matrix for $\phi$ is the
same under models (\ref{eq.model.sansperiode}) and (\ref
{eq.modele.avecperiode}),
that is,
\[
\widetilde C_d[\phi]=C_d[\phi]. %
\]
In particular, this is true if $d$ is a symmetric design.
\end{proposition}

\begin{pf}
The method of proof of Lemma~\ref{lemma.Lstarsymmetric} shows that the
matrix $\widetilde L$ used for minimising may be chosen to satisfy
$P_{\bar\sigma}' \widetilde L P_\sigma= \widetilde L$
for all $\sigma$ in $G_d$. This means that $L=L_\sigma$ and
$N = N_\sigma=NP_\sigma$ for all $\sigma$ in $G_d$. If $NP_\sigma= N$
for all $\sigma$ in $G_d$, and $G_d$ is transitive, then every row of
$N$ is a multiple of $\mathbb{I}_t'$.

We have $C_{d12} = X_d'\omega_B^\perp A = X_d' A Q_k$.
Lemma~\ref{lemma:strong} shows that if $L=L_\sigma$ then
$L'C_{d12} = L_\sigma'X_d'A Q_k =
L'_\sigma P_{\widetilde{\sigma}}'X_d' A Q_k =
P_\sigma' L' C_{d12}$.
If $G_d$ is transitive, then every column of 
$L'C_{d12}$ is a multiple of $\mathbb{I}_t$.

Therefore, the expression in (\ref{eq.CminLperiode}) is equal to
$L'C_d[\xi]L + c(L,N)J_t$ for some scalar $c(L,N)$. Hence
\begin{eqnarray*}
\widetilde C_d[\phi]= Q_t \widetilde C_d[
\phi] Q_t &=&Q_t \Bigl( \min_{(L'\mid N')'\in\lsett}
L'C_d[\xi]L +c(L,N)J_t \Bigr)
Q_t
\\
& = & \min_{(L' \mid N')'\in\lsett} 
 \bigl( Q_tL'C_d[
\xi]LQ_t \bigr)
\\
& = & Q_t \Bigl(\min_{L\in\lset}L'C_d[
\xi]L \Bigr) Q_t
\\
& = & Q_t C_d[\phi] Q_t =
C_d[\phi].
\end{eqnarray*}
\upqed\end{pf}

For any design $d$ whose $G_d$ is doubly transitive,
$C_d[\phi]$ is completely symmetric (replace $S_t$ by $G_d$ in the
proof of
Lemma~\ref{lemma.cphics}).
Double transitivity implies strong balance on the periods, so
then $\widetilde C_d[\phi]$ is also completely symmetric, by
Proposition~\ref{prop.equalinfomatperiode}.
In Section~\ref{subsec:small} we give some examples that show that strong
balance on the periods is not sufficient for $\widetilde C_d[\phi]$ to be
completely symmetric.

The results obtained in this section also hold for approximate designs.
Since the restriction of $A$ to a single sequence is equal to $I_k$,
for an exact designs $d$ we have
\[
\widetilde C_d[ \theta]= n\sum_{s\in\mathcal S}
\pi_d(s)\pmatrix{
X'_{s} Q_k X_s &
X'_s Q_k
\vspace*{2pt}\cr
Q_k X_s & Q_k 
}.
\]
%
This expression can also be used for
approximate designs. Moreover, in the definition of a design being
strongly balanced on the periods ``equally often'' may be replaced by ``in
the same proportions'' and ``number of times'' by ``proportion of times.''
Then the proofs of Lemma~\ref{lemma:strong} and
Proposition~\ref{prop.equalinfomatperiode} can be easily adapted to
approximate designs by replacing $A' X_d$ by $n\sum_{s} \pi_d(s)X_{s}$,
replacing $X_d'\omega_B^\perp A$ by $\sum_{s} \pi_d(s)X_{s}Q_k$, and
so on.

\section{Universally optimal approximate designs}
\label{section.universal}

From \citet{kief:1975}, a design $d^*$ for which the information matrix
$C_{d^*}[\phi]$ is completely symmetric and that maximises the trace
of $C_d[\phi]$ over all the designs $d$
for $t$ treatments using $n$ subjects for $k$ periods
is universally optimal.

\subsection{Condition for optimal designs}
The following proposition shows that a universally optimal approximate
design may be sought among symmetric designs.

\begin{proposition}
A symmetric design for which the trace of the information matrix is
maximal among the class of symmetric designs is universally optimal
among all possible approximate designs.
\end{proposition}

\begin{pf}
For any design $d$, taking the trace in (\ref{eq.Cphidbar}), we have
$\operatorname{tr}(C_{\bar d}[\phi])\geq\operatorname{tr}(C_{d}[\phi
])$. Since, by Lemma~\ref{lemma.cphics}, $C_{\bar d}[\phi]$ is
completely symmetric, $\bar d$ is always better than $d$ with respect
to universal optimality. If $d^*$ maximises the trace among the set of
symmetric designs, then for any design $d$, $\operatorname{tr}(C_{{d^*}}[\phi])\geq\operatorname{tr}(C_{\bar d}[\phi])\geq
\operatorname{tr}(C_{d}[\phi])$. Since $
C_{{d^*}}[\phi]$ is completely symmetric and maximises the trace, $d^*$
is universally optimal.
\end{pf}

For any sequence $s$, and $1\leq p,q \leq7$, put
$c_{spq}=\operatorname{tr} (L_{(p)}' C_s[\xi] L_{(q)} )$.
Then combining (\ref{eq.cphiLstar}), (\ref{eq.infomatphi}), and
(\ref{eq.Lstarreduc}), we have for a symmetric design,
\[
\operatorname{tr}\bigl(C_{ d}[\phi]\bigr)=\min_{\gamma_2,\ldots,\gamma_6}
\sum_{s\in\mathcal S}n \pi_d(s)\sum
_{p=1}^6\sum_{q=1}^6
\gamma_p\gamma_q c_{spq}\qquad \mbox{with }
\gamma_1=1. %
\]

\begin{lemma}
\label{lemme.csijsigma}
For a sequence $s$ and a permutation $\sigma$ on the treatment labels,
we have
\[
c_{s_\sigma pq}=c_{spq}. %
\]
\end{lemma}

\begin{pf}
\begin{eqnarray*}
 c_{s_\sigma pq} &=& \operatorname{tr}
\bigl(P_\sigma' L_{(p)}'
C_{s_\sigma}[\xi] L_{(q)}P_\sigma\bigr)\qquad  \mbox{since }
\operatorname{tr}(AB)=\operatorname{tr}(BA),
\\
&=&\operatorname{tr}\bigl(P_\sigma' L_{(p)}'P_{\widetilde\sigma}
C_s[\xi] P'_{\widetilde\sigma} L_{(q)}P_\sigma
\bigr) \qquad \mbox{by }(\ref{eq.Cxidsigma}),
\\
&=& \operatorname{tr}\bigl(L_{(p)}' C_s[\xi]
L_{(q)}\bigr)=c_{spq} \qquad \mbox{by }(\ref{eq.invarianceLi}).
\end{eqnarray*}
\upqed\end{pf}

Two sequences are said to be \emph{equivalent} if one can be obtained
from the other one by some permutation of treatment labels. We denote by
$\ecset$ the set of all possible \emph{equivalence classes}.
From Lemma~\ref{lemme.csijsigma}, $c_{spq}$ depends only on the equivalence
class $\ell$ to which $s$ belongs, and will be therefore denoted
$c_{\ell pq}$.
To each equivalence class~$\ell$, we may also associate the
nonnegative convex
quadratic polynomial with five variables $\gamma=(\gamma_2,\ldots,\gamma_6)$,
\[
h_\ell(\gamma)=\sum_{p=1}^6
\sum_{q=1}^6 \gamma_p
\gamma_q c_{\ell pq} \qquad\mbox{where }\gamma_1=1.
\]
For a symmetric design, we may write $\pi_\ell$ for the proportion of
sequences which are in the equivalence class $\ell$. Then
\[
\operatorname{tr}\bigl(C_{ d}[\phi]\bigr)=\min_{\gamma}
\sum_{\ell\in\ecset} n \pi _\ell h_\ell(
\gamma). %
\]
Therefore, we have the following proposition:

\begin{proposition}
An approximate symmetric design $d^*$ with proportions $\{\pi^*_\ell\}
_{\ell\in\ecset}$ that achieves
%
\begin{equation}
\label{eq.maxmin} \max_{\{\pi_\ell\}_{\ell\in\ecset}} \min_\gamma\sum
_{\ell\in\ecset} \pi _\ell h_\ell(\gamma)
\end{equation}
is universally optimal for $\phi$ among all possible designs.
\end{proposition}

\subsection{Determination of optimal proportions}
\label{subsection.procedure}
Each equivalence class of sequences is defined by a partition of the set
$\{1, 3, \ldots, k\}$ into at most $t$ parts. If $t\geq k$, the number
of such partitions is the\vadjust{\goodbreak} Bell number $B_k$, which grows with $k$ more
than exponentially [\citet{cam:1994}, Chapter~3]. Thus it is not realistic
to solve the maximin problem in (\ref{eq.maxmin}) by hand.

It seems intuitive that sequences in an optimal symmetric design should
satisfy two contradictory conditions: for accurate estimation of total
effects, each treatment should be preceded by itself a large number of
times; while, for efficiency in allowing for subjects, the replications
within each sequence should be as equal as possible. As a compromise,
this suggests sequences in which all occurrences of each treatment are
in a run of consecutive periods. Indeed, in our numerical results in
Section~\ref{section.results}, all seqeuences in the optimal designs
have this form. Each equivalence class of such sequences is defined by
a so-called \emph{composition} of $k$. However, the number of
compositions of $k$ is $2^{k-1}$
[\citet{cam:1994}, Chapter~4], so, even if we restrict ourselves to such
sequences, a hand search is still not realistic.

We propose now the following method derived from \citet{kush:1997}.
Consider
\[
h^*(\gamma)=\max_{\ell\in\ecset}h_\ell(\gamma).
\]
We use the following procedure.
\begin{longlist}[\textit{Step} 1.]
\item[\textit{Step} 1.] Find $\gamma^*$ that minimises the function $h^*(\gamma
)$, and denote $h^*=h^*(\gamma^*)$ the minimum.
\item[\textit{Step} 2.] Select the classes $\ell$ of sequences such that
$h_\ell(\gamma^*)=h^*$, and denote $\ecset^*$ this set.
\item[\textit{Step} 3.] Solve in $\{\pi_\ell\mid\ell\in\ecset^*\}$ the linear system,
$\sum_{\ell\in\ecset^*}\pi_\ell\frac{d h_\ell}{d\gamma} (\gamma^*) =0$,
for $0<\pi_\ell<1$ and $\sum_{\ell\in\ecset} \pi_\ell=1$; denote
$\pi^*=\{\pi^*_\ell\mid\ell\in\ecset^*\}$ the solution
(not necessarily unique).
\item[\textit{Step} 4.] Give the symmetric designs such that
$\pi_\ell=\pi^*_\ell$ for $\ell\in\ecset^*$ and $\pi_\ell=0$
otherwise; these designs are universally optimal.
\end{longlist}

Step 1 is the most challenging. However,
since $h^*(\gamma)$
is a convex function, any standard optimisation algorithm gives
accurate values for $\gamma^*$ and $h^*$ in a short time, even if the
number of possible classes is large. When supported by the software, we
used an exact optimisation algorithm to obtain the values of $\gamma^*$.

For step 2, the optimal sequences are part of the information found in
step 1. Since $\mathcal{C}^*$ is usually rather small, step 3 simply
involves inverting a small square matrix whose entries have been found
in step 1. Step 4 then reports the results.

\section{Examples of optimal and efficient designs}
\label{section.results}

For some values of $k$ and $t$, we give optimal approximate designs for
$\phi$.
For each given $k$, the first table gives the optimal proportions, and
the second table gives the efficiency factor for a symmetric design
generated by a single sequence.

Consider a real-valued criterion $\psi(C_d[\phi])$ which is concave,
nondecreasing in $C_d[\phi]$ with respect to the Loewner ordering,\vadjust{\goodbreak} and
invariant under simultaneous permutations of rows and columns.
From \citet{kief:1975},
there is an approximate design $d^*$ which maximises $\psi(C_d[\phi])$
over the set of approximate designs with the same values of $k$ and $t$.
The efficiency factor of a design $d$ for criterion $\psi$
can therefore be defined by
\[
\mathit{eff}_\psi(d)=\frac{\psi(C_d[\phi])}{\psi(C_{d^*}[\phi])}. %
\]
%
For $\psi(C)=\operatorname{tr}(C)$, we simply write
%
\begin{equation}
\label{eq.defeff} \mathit{eff} (d)=\frac{\operatorname{tr}(C_d[\phi
])}{\operatorname{tr}(C_{d^*}[\phi])}.
\end{equation}
When $C_d[\phi]$ is completely symmetric,
$\mathit{eff} (d)$ is also the efficiency factor for the well-known
$D$-, $A$- and $E$-criteria; see \citet{shah:sinh:1989} or \citet{drui:2004}.

In our tables,
we write $0^+$ or $1^-$ when a value is 
within $0.005$ of $0$, $1$, respectively. For some values of $k$ and $t$ the
optimal proportions have been calculated with formal calculus when tractable;
all others have been obtained by numerical optimisation.

The values $h^*$ displayed correspond to those defined
in Section~\ref{subsection.procedure} for an optimal design.
The information matrix for a symmetric optimal approximate design with
$n$ subjects is therefore
\[
C_d[\phi]= \frac{n h^*}{t-1} Q_t. 
\]

\subsection{3 periods}
Optimal proportions for some values of $t$:\vspace*{12pt}

\noindent\tabcolsep=0pt
\begin{tabular*}{\textwidth}{@{\extracolsep{\fill}}lccccccccccccccc@{}}
\hline
$\bolds{t}$ &\textbf{2}&\textbf{3}&\textbf{4}&\textbf{5}&\textbf{6}&\textbf{7}&\textbf{8}&\textbf{9}&\textbf{10}&\textbf{11}&\textbf{12}&\textbf{13}&\textbf{14}&\textbf{15}&\textbf{16}\\
\hline
Prop. [ 1 1 2 ] &$\tfrac{1}{2}$ &
$\tfrac{5}{13}$&$\tfrac{1}{3}$& $\tfrac
{7}{23}$ &$\tfrac{2}{7}$&
$\tfrac{3}{11}$& $\tfrac{5}{19}$&$\tfrac{11}{43}$&$\tfrac
{1}{4}$&
$\tfrac{13}{53}$&$\tfrac{7}{29}$&$\tfrac{5}{21}$&$\tfrac{4}{17}$&
$\tfrac
{17}{73}$&$\tfrac{3}{13}$
\\[3pt]
Prop. [ 1 2 2 ]& $\tfrac{1}{2}$ &
$\tfrac{8}{13}$& $\tfrac{2}{3}$ & $\tfrac{16}{23}$ &$\tfrac{5}{7}$&
$\tfrac{8}{11}$& $\tfrac
{14}{19}$&$\tfrac{32}{43}$&$\tfrac{3}{4}$&
$\tfrac{40}{53}$&$\tfrac{22}{29}$&$\tfrac
{16}{21}$&$\tfrac{13}{17}$&
$\tfrac{56} {73}$&$\tfrac{10}{13}$
\\[3pt]
$h^*$& $\tfrac{1}{3}$&$\tfrac{16}{39}$&
$\tfrac{4}{9}$&$\tfrac{32}{69}$& $\tfrac
{10}{21}$&$\tfrac{16}{33}$&
$\tfrac{28}{57}$&$\tfrac{64}{129}$&$\tfrac{1}{2}$& $\tfrac{80}{159}$&
$\tfrac{44}{87}$&$\tfrac{32}{63}$&$\tfrac{26}{51}$&$\tfrac
{112}{219}$&
$\tfrac{20}{39}$
\\
\hline
\end{tabular*}\vspace*{12pt}

Efficiency of symmetric designs generated by a single sequence:\vspace*{12pt}

\noindent\tabcolsep=0pt
\begin{tabular*}{\textwidth}{@{\extracolsep{\fill}}lccccccccccccccc@{}}
\hline
$\bolds{t}$ &\textbf{2}&\textbf{3}&\textbf{4}&\textbf{5}&\textbf{6}&\textbf{7}&\textbf{8}&\textbf{9}&\textbf{10}&\textbf{11}&\textbf{12}&\textbf{13}&\textbf{14}&\textbf{15}&\textbf{16}\\
\hline
{Eff. } [ 1 1 2 ] &0& 0&0&0&0&0&0&0&0&0&0&0&0&0&0
\\
{Eff. } [ 1 2 2 ]&0&0.61&0.75&0.81&0.84&0.86&0.87&0.88&0.89&0.89&0.90&0.90&0.91&0.91&0.91
\\
\hline
\end{tabular*}\vspace*{12pt}

 Example of universally optimal design for $t=4$:
\[
\lleft( %
\begin{array} {cccccccccccccccccccccccccccccccccccc} 1 & 1 &
1 & 1 & 1 & 1 & 2 & 2 & 2 & 2 & 2 & 2 & 3 & 3 & 3 & 3 & 3 & 3 & 4 & 4 & 4 & 4 & 4
& 4 & 1 & 1 & 1 & 2 & 2 & 2 & 3 & 3 & 3 & 4 & 4 & 4
\\
2 & 2 & 3 & 3 & 4 & 4 & 1 & 1 & 3 & 3 & 4 & 4 & 1 & 1 & 2 & 2 & 4 & 4 & 1 & 1 &
2 & 2 & 3 & 3 & 1 & 1 & 1 & 2 & 2 & 2 & 3 & 3 & 3 & 4 & 4 & 4
\\
2 & 2 & 3 & 3 & 4 & 4 & 1 & 1 & 3 & 3 & 4 & 4 & 1 & 1 & 2 & 2 & 4 & 4 & 1 & 1 &
2 & 2 & 3 & 3 & 2 & 3 & 4 & 1 & 3 & 4 & 1 & 2 & 4 & 1 & 2 & 3 \end{array}
 \rright) %
\]
%

\subsection{4 periods}
The optimal approximate designs are generated by the single sequence $[1\ 1\ 2\ 2]$ for $2\leq t\leq30$. It is conjectured that this is true
for any value of~$t$.
\subsection{5 periods}

Optimal proportions for some values of $t$:\vspace*{12pt}

%
\noindent\tabcolsep=0pt
\begin{tabular*}{\textwidth}{@{\extracolsep{\fill}}lcccccccccccc@{}}
\hline
$\bolds{t}$ &\textbf{2}&\textbf{3}&\textbf{4}&\textbf{5}&\textbf{6}&\textbf{7}&\textbf{8}&\textbf{9}&\textbf{10}&\textbf{15}&\textbf{20}&\textbf{30}\\
\hline
Prop. [ 1 1 2 2 2 ]&$\tfrac{1}{2}$ &
$\tfrac{7}{9}$& $\tfrac{17}{19}$ & $\tfrac{47}{49}$ &0.98&0.98&0.98&0.98&0.98&0.97&0.97&0.97
\\[3pt]
Prop. [ 1 1 1 2 2 ] &$\tfrac{1}{2}$&
$\tfrac{2}{9}$&$\tfrac{2}{19}$ &$\tfrac{2}{49}$ &0\phantom{00.}&0\phantom{00.}&0\phantom{00.}&0\phantom{00.}&0\phantom{00.}&0\phantom{00.}&0\phantom{00.}&0\phantom{00.}
\\[3pt]
Prop. [ 1 1 2 3 3 ]&0&0&0&0&0.02&0.02&0.02&0.02&0.02&0.03&0.03&0.03
\\[3pt]
$h^*$& $\tfrac{7}{5}$& $\tfrac{68}{45}$ &
$\tfrac{148}{95}$ & $\tfrac{388}{245}$ &1.60&1.61&1.62&1.63&1.63&1.64&1.65&1.66
\\
\hline
\end{tabular*}\vspace*{12pt}

 Efficiency of symmetric designs generated by a single sequence:\vspace*{12pt}

\noindent\tabcolsep=0pt
\begin{tabular*}{\textwidth}{@{\extracolsep{\fill}}lcccccccccccc@{}}
\hline
$\bolds{t}$ &\textbf{2}&\textbf{3}&\textbf{4}&\textbf{5}&\textbf{6}&\textbf{7}&\textbf{8}&\textbf{9}&\textbf{10}&\textbf{15}&\textbf{20}&\textbf{30}\\
\hline
Eff. [ 1 1 2 2 2 ] &0.95&0.99&0.998 &
$1^-$\phantom{0}&$1^-$\phantom{0}&$1^-$\phantom{0}&$1^-$\phantom{0}&$1^-$\phantom{0}&$1^-$\phantom{0}&$1^-$\phantom{0}&$1^-$\phantom{0}&$1^-$\phantom{0}
\\
Eff. [ 1 1 1 2 2 ]&0.95 &0.91&0.89 &
0.88&0.87&0.87&0.87&0.87&0.87& 0.86&0.86&0.85
\\
Eff. [ 1 1 2 3 3 ]& --&0.77&0.82&0.84&0.85&0.86&0.86&
0.86& 0.86& 0.87&0.88&0.88
\\
\hline
\end{tabular*}\vspace*{6pt}\

 Example of universally optimal symmetric design for $t=3$:

\begin{eqnarray*}
\left(
\begin{array}{llllllllllllllllllllllllllll}
1 & 1 & 1 & 1 & 1 & 1 & 1 & 1 & 1 & 1 & 1 & 1 & 1 & 1 & 2 & 2 & 2 & 2
& 2 & 2 & 2 & 2 & 2 & 2 & 2 & 2 & 2 & 2 \\
1 & 1 & 1 & 1 & 1 & 1 & 1 & 1 & 1 & 1 & 1 & 1 & 1 & 1 & 2 & 2 & 2 & 2
& 2 & 2 & 2 & 2 & 2 & 2 & 2 & 2 & 2 & 2 \\
2 & 2 & 2 & 2 & 2 & 2 & 2 & 3 & 3 & 3 & 3 & 3 & 3 & 3 & 1 & 1 & 1 & 1
& 1 & 1 & 1 & 3 & 3 & 3 & 3 & 3 & 3 & 3 \\
2 & 2 & 2 & 2 & 2 & 2 & 2 & 3 & 3 & 3 & 3 & 3 & 3 & 3 & 1 & 1 & 1 & 1
& 1 & 1 & 1 & 3 & 3 & 3 & 3 & 3 & 3 & 3 \\
2 & 2 & 2 & 2 & 2 & 2 & 2 & 3 & 3 & 3 & 3 & 3 & 3 & 3 & 1 & 1 & 1 & 1
& 1 & 1 & 1 & 3 & 3 & 3 & 3 & 3 & 3 & 3
\end{array}
\right.
 \\
\left.
\begin{array}{llllllllllllllllllllllllll}
3 & 3 & 3 & 3 & 3 & 3 & 3 & 3 & 3 & 3 & 3 & 3 & 3 & 3 & 1 & 1 & 1 & 1
& 2 & 2 & 2 & 2 & 3 & 3 & 3 & 3 \\
3 & 3 & 3 & 3 & 3 & 3 & 3 & 3 & 3 & 3 & 3 & 3 & 3 & 3 & 1 & 1 & 1 & 1
& 2 & 2 & 2 & 2 & 3 & 3 & 3 & 3 \\
1 & 1 & 1 & 1 & 1 & 1 & 1 & 2 & 2 & 2 & 2 & 2 & 2 & 2 & 1 & 1 & 1 & 1
& 2 & 2 & 2 & 2 & 3 & 3 & 3 & 3 \\
1 & 1 & 1 & 1 & 1 & 1 & 1 & 2 & 2 & 2 & 2 & 2 & 2 & 2 & 2 & 2 & 3 & 3
& 1 & 1 & 3 & 3 & 1 & 1 & 2 & 2 \\
1 & 1 & 1 & 1 & 1 & 1 & 1 & 2 & 2 & 2 & 2 & 2 & 2 & 2 & 2 & 2 & 3 & 3
& 1 & 1 & 3 & 3 & 1 & 1 & 2 & 2
\end{array}
\right)\hspace*{-5pt}
\end{eqnarray*}

\subsection{6 periods}

Optimal proportions for some values of $t$:\vspace*{9pt}

\noindent\tabcolsep=0pt
\begin{tabular*}{\textwidth}{@{\extracolsep{\fill}}lcccccccccccc@{}}
\hline
$\bolds{t}$ &\textbf{2}&\textbf{3}&\textbf{4}&\textbf{5}&\textbf{6}&\textbf{7}&\textbf{8}&\textbf{9}&\textbf{10}&\textbf{15}&\textbf{20}&\textbf{30}\\
\hline
Prop. [ 1 1 1 2 2 2 ] & 1&0.81&0.66 &0.55
&0.48&0.42&0.38&0.35&0.32&0.23&0.19&0.15
\\
Prop. [ 1 1 2 2 3 3 ]&0& 0.19&0.34 & 0.45
&0.52&0.58 &0.62&0.65&0.68&0.77&0.81&0.85
\\
$h^*$& 2 & 2.11&2.16 & 2.19 &2.21&2.22&2.23&2.24&2.25&2.26&2.27&2.28
\\
\hline
\end{tabular*}\vspace*{9pt}\vadjust{\goodbreak}

 Efficiency of symmetric designs generated by a single sequence:\vspace*{9pt}

\noindent\tabcolsep=0pt
\begin{tabular*}{\textwidth}{@{\extracolsep{\fill}}lcccccccccccc@{}}
\hline
$\bolds{t}$ &\textbf{2}&\textbf{3}&\textbf{4}&\textbf{5}&\textbf{6}&\textbf{7}&\textbf{8}&\textbf{9}&\textbf{10}&\textbf{15}&\textbf{20}&\textbf{30}\\
\hline
Eff. [ 1 1 1 2 2 2 ] & 1&0.99&0.99&0.98&0.98&0.97&0.97&0.97&0.97&0.97&0.96&0.96
\\
Eff. [ 1 1 2 2 3 3 ]& --& 0.95&0.97&0.98&0.99&0.99&0.99&0.99&$1^-$\phantom{0}&$1^-$\phantom{0}&$1^-$\phantom{0}&$1^-$\phantom{0}
\\
\hline
\end{tabular*}\vspace*{9pt}

\subsection{7 periods}

Optimal proportions for some values of $t$:\vspace*{9pt}

\noindent\tabcolsep=0pt
\begin{tabular*}{\textwidth}{@{\extracolsep{\fill}}lccccc@{}}
\hline
$\bolds{t}$ & \textbf{3}&\textbf{4}&\textbf{5}&\textbf{6}&\multicolumn{1}{c@{}}{$\bolds{7\leq t\leq30}$}
\\
\hline
Prop. [ 1 1 1 2 2 2 2 ] & 0.57&0.19&0\phantom{00.}&0\phantom{00.}&0\phantom{00.}
\\
Prop. [ 1 1 1 2 2 3 3 ] &0\phantom{00.}&0\phantom{00.}&0.09&$0^+$\phantom{0}&0\phantom{00.}
\\
Prop. [ 1 1 2 2 3 3 3 ]& 0.43&0.81&0.91&$1^-$\phantom{0}&1\phantom{00.}
\\
$h^*$& 2.60&2.70&2.76&2.80&2.82
\\
\hline
\end{tabular*}\vspace*{9pt}

 Efficiency of symmetric designs generated by a single sequence:\vspace*{9pt}

\noindent\tabcolsep=0pt
\begin{tabular*}{\textwidth}{@{\extracolsep{\fill}}lccccc@{}}
\hline
$\bolds{t}$ & \textbf{3}&\textbf{4}&\textbf{5}&\textbf{6}&\multicolumn{1}{c@{}}{\textbf{7}}
\\
\hline
Eff.  [ 1 1 1 2 2 2 2 ] & 0.98
&0.96&0.95&0.94&0.94
\\
Eff. [ 1 1 1 2 2 3 3 ] &0.98&0.99&0.98&0.98&0.98
\\
Eff. [ 1 1 2 2 3 3 3 ]& 0.98& $1^-$\phantom{0}&$1^-$\phantom{0}&$1^-$\phantom{0}&1\phantom{00.}
\\
\hline
\end{tabular*}

\subsection{Efficient designs with $t(t-1)$ subjects}
\label{subsec:small}
For $k=6$ or $k=7$, we saw that efficient symmetric designs may be
obtained from single sequences having
three treatments by permuting all the treatment labels.
Such designs require $t(t-1)(t-2)$ subjects, which may be too large.
We can construct efficient designs that are strongly balanced on the
periods, are generated by a single sequence, and require only $t(t-1)$
subjects, as follows.
%
\begin{enumerate}[\textit{Step} 1.]
\item[\textit{Step} 1.] We start from a balanced incomplete-block design with
block-size $3$ and $t$ treatments such that for any two different
periods $j_1$ and $j_2$ and any two different treatments $u$ and $v$,
there exists exactly one subject that receives treatment $u$ in period
$j_1$ and treatment $v$ in period
$j_2$. [This is called an orthogonal array of type I and strength two;
see \citet{rao:1961}.]
\begin{itemize}
\item If $t$ is odd, use all the triplets
$[u,u+v, u+2v]$ modulo $t$, for $u=0, \ldots, t-1$ and $v=1, \ldots, t-1$.
%
\item If $t$ is even, use the preceding construction for $t-1$ and
replace each triplet of the form
$[u,u+1,u+2]$
by the three sequences
$[t,u+1,u+2]$, $[u,t,u+2]$ and $[u,u+1,t]$.
\end{itemize}
\item[\textit{Step} 2.] Then we construct a design with $k$ periods by
replicating the three treatments in each triplet in such a way that we
obtain a sequence
in the same equivalence class as the one that
generates the efficient design.
\end{enumerate}
For example, take $k=7$ and $t=5$ with generating sequence $[ 1\ 1\ 2\ 2\ 3\ 3\ 3 ]$.

The starting design with three periods is

\[
\lleft( %
\begin{array} {cccccccccccccccccccc} 1 & 1 & 1 & 1 & 2 & 2 &
2 & 2 & 3 & 3 & 3 & 3 & 4 & 4 & 4 & 4 & 5 & 5 & 5 & 5
\\
2 & 3 & 4 & 5 & 3 & 4 & 5 & 1 & 4 & 5 & 1 & 2 & 5 & 1 & 2 & 3 & 1 & 2 & 3 & 4
\\
3 & 5 & 2 & 4 & 4 & 1 & 3 & 5 & 5 & 2 & 4 & 1 & 1 & 3 & 5 & 2 & 2 & 4 & 1 & 3
\end{array} %
 \rright). %
\]

The resulting design with seven periods generated by $[1\ 1\ 2\ 2\ 3\ 3\ 3]$ is
%
\[
\lleft( %
\begin{array} {cccccccccccccccccccc} 1 & 1 & 1 & 1 & 2 & 2 &
2 & 2 & 3 & 3 & 3 & 3 & 4 & 4 & 4 & 4 & 5 & 5 & 5 & 5
\\
1 & 1 & 1 & 1 & 2 & 2 & 2 & 2 & 3 & 3 & 3 & 3 & 4 & 4 & 4 & 4 & 5 & 5 & 5 & 5
\\
2 & 3 & 4 & 5 & 3 & 4 & 5 & 1 & 4 & 5 & 1 & 2 & 5 & 1 & 2 & 3 & 1 & 2 & 3 & 4
\\
2 & 3 & 4 & 5 & 3 & 4 & 5 & 1 & 4 & 5 & 1 & 2 & 5 & 1 & 2 & 3 & 1 & 2 & 3 & 4
\\
3 & 5 & 2 & 4 & 4 & 1 & 3 & 5 & 5 & 2 & 4 & 1 & 1 & 3 & 5 & 2 & 2 & 4 & 1 & 3
\\
3 & 5 & 2 & 4 & 4 & 1 & 3 & 5 & 5 & 2 & 4 & 1 & 1 & 3 & 5 & 2 & 2 & 4 & 1 & 3
\\
3 & 5 & 2 & 4 & 4 & 1 & 3 & 5 & 5 & 2 & 4 & 1 & 1 & 3 & 5 & 2 & 2 & 4 & 1 & 3
\end{array} %
 \rright). %
\]
%

The following table displays the A-, D-, E-efficiency factors for
designs with $6$ periods and $t(t-1)$ subjects generated by the
sequence $[ 1\ 1\ 2\ 2\ 3\ 3 ]$ using
the method described above. The efficiency factors are 
given relative to universally optimal approximate designs.\vspace*{6pt}

%
\noindent\tabcolsep=0pt
\begin{tabular*}{\textwidth}{@{\extracolsep{\fill}}lccccccc@{}}
\hline
$\bolds{t}$ & \textbf{4}&\textbf{5}&\textbf{6}&\textbf{7} &\textbf{8}&\textbf{9}&\textbf{10}
\\
\hline
A-efficiency & 0.951 & 0.977 & 0.973 & 0.978 &
0.974 & 0.970&0.968
\\
D-efficiency & 0.951 & 0.977 & 0.973 & 0.978 &
0.974&0.970&0.968
\\
E-efficiency & 0.951 & 0.977 & 0.951 & 0.978 &
0.950 &0.950&0.949
\\
\hline
\end{tabular*}\vspace*{6pt}

We may note that this method is interesting only for $t=7$ or $t=8$.
For the other values of $t$, the symmetric design with $t(t-1)$
subjects generated by the sequence $[ 1\ 1\ 1\ 2\ 2\ 2 ]$ is more efficient.

The following table displays the A-, D-, E-efficiency factors for
designs with $7$ periods and $t(t-1)$ subjects generated by the
sequence $[ 1\ 1\ 2\ 2\ 3\ 3\ 3 ]$ using
the method described above.\vspace*{6pt}

\noindent\tabcolsep=0pt
\begin{tabular*}{\textwidth}{@{\extracolsep{\fill}}lccccccc@{}}
\hline
$\bolds{t}$ & \textbf{4}&\textbf{5}&\textbf{6}&\textbf{7} &\textbf{8}&\textbf{9}&\textbf{10}
\\
\hline
A-efficiency & 0.974 &0.990& 0.982 & 0.983 & 0.978
& 0.973 &0.971
\\
D-efficiency & 0.974 &0.990& 0.982& 0.983 & 0.978
& 0.973 & 0.971
\\
E-efficiency & 0.974 &0.990& 0.961& 0.983 & 0.955
& 0.954 & 0.954
\\
\hline
\end{tabular*}\vspace*{6pt}

For $t=4,5,7$, the information matrices are completely symmetric. For
$t\geq4$ and when the number of subjects is $t(t-1)$, these designs
are preferable to symmetric designs generated by the sequence $[ 1\ 1\ 1\ 2\ 2\ 2\ 2 ]$.

If $t=4$ or $t$ is an odd prime, this method always gives a design $d$
for which $G_d$ is doubly transitive, and so $\widetilde C_d[\phi]$ is
completely symmetric. If $t$ is any prime power, there is a second
method which gives a design $d$ in $t(t-1)$ periods for which $G_d$ is
completely symmetric.

\begin{longlist}[\textit{Step} 1.]
\item[\textit{Step} 1.] Identify the treatments with the elements of the finite field
$\mathrm{GF}(t)$ of order $t$.
\item[\textit{Step} 2.] Form any triplet $[x,y,z]$ of distinct treatments.
\item[\textit{Step} 3.] Use this to produce all triplets of the form
$[ax+b,ay+b,az+b]$ for which $a$ and $b$ are in $\mathrm
{GF}(t)$ and $a\ne0$.
\item[\textit{Step} 4.] Use these triplets to construct a design from the desired
sequence just as in the previous method.
\end{longlist}

For example, when $t=8$, one correspondence between $\{1, \ldots, 8\}$ and
$\mathrm{GF}(8)$ gives the following starting design with
three periods:
\begin{eqnarray*}
\left(
\begin{array} {cccccccccccccccccccccccccccccccc} 8 & 7 & 1 & 3 & 2 & 6 & 4 & 5 & 8
& 1 & 2 & 4 & 3 & 7 & 5 & 6 & 8 & 2 & 3 & 5 & 4 & 1 & 6 & 7 & 8 & 3 & 4 & 6 & 5 &
2 & 7 & 1
\\
7 & 8 & 3 & 1 & 6 & 2 & 5 & 4 & 1 & 8 & 4 & 2 & 7 & 3 & 6 & 5 & 2 & 8 & 5 & 3 &
1 & 4 & 7 & 6 & 3 & 8 & 6 & 4 & 2 & 5 & 1 & 7
\\
1 & 3 & 8 & 7 & 4 & 5 & 2 & 6 & 2 & 4 & 8 & 1 & 5 & 6 & 3 & 7 & 3 & 5 & 8 & 2 &
6 & 7 & 4 & 1 & 4 & 6 & 8 & 3 & 7 & 1 & 5 & 2 \end{array} %
\right.
\\
\left.
\begin{array} {cccccccccccccccccccccccc} 8 & 4 & 5 & 7 & 6
& 3 & 1 & 2 & 8 & 5 & 6 & 1 & 7 & 4 & 2 & 3 & 8 & 6 & 7 & 2 & 1 & 5 & 3 & 4
\\
4 & 8 & 7 & 5 & 3 & 6 & 2 & 1 & 5 & 8 & 1 & 6 & 4 & 7 & 3 & 2 & 6 & 8 & 2 & 7 &
5 & 1 & 4 & 3
\\
5 & 7 & 8 & 4 & 1 & 2 & 6 & 3 & 6 & 1 & 8 & 5 & 2 & 3 & 7 & 4 & 7 & 2 & 8 & 6 &
3 & 4 & 1 & 5 \end{array}\right).\hspace*{-10pt}
\end{eqnarray*}

The design obtained from this starting design and the generating sequence
$[ 1\ 1\ 2\ 2\ 3\ 3 ]$, respectively, $[ 1\ 1\ 2\ 2\ 3\ 3\ 3 ]$,
has efficiency factor equal to $0.977$, respectively, to $0.981$.

For $t=9$, we obtain the following starting design:
\begin{eqnarray*}
\left(
\begin{array} {cccccccccccccccccccccccccccccccccccc} 1 & 1 & 2 & 2 & 3 & 3 & 4 & 4
& 5 & 5 & 6 & 6 & 7 & 7 & 8 & 8 & 9 & 9 & 1 & 1 & 4 & 4 & 7 & 7 & 2 & 2 & 5 & 5 &
8 & 8 & 3 & 3 & 6 & 6 & 9 & 9
\\
2 & 3 & 1 & 3 & 1 & 2 & 5 & 6 & 4 & 6 & 4 & 5 & 8 & 9 & 7 & 9 & 7 & 8 & 4 & 7 &
1 & 7 & 1 & 4 & 5 & 8 & 2 & 8 & 2 & 5 & 6 & 9 & 3 & 9 & 3 & 6
\\
3 & 2 & 3 & 1 & 2 & 1 & 6 & 5 & 6 & 4 & 5 & 4 & 9 & 8 & 9 & 7 & 8 & 7 & 7 & 4 &
7 & 1 & 4 & 1 & 8 & 5 & 8 & 2 & 5 & 2 & 9 & 6 & 9 & 3 & 6 & 3 \end{array}
\right.
\\
\left.
\begin{array} {cccccccccccccccccccccccccccccccccccc} 1 & 1
& 5 & 5 & 9 & 9 & 2 & 2 & 6 & 6 & 7 & 7 & 3 & 3 & 4 & 4 & 8 & 8 & 1 & 1 & 6 & 6 &
8 & 8 & 2 & 2 & 4 & 4 & 9 & 9 & 3 & 3 & 5 & 5 & 7 & 7
\\
5 & 9 & 1 & 9 & 1 & 5 & 6 & 7 & 2 & 7 & 2 & 6 & 4 & 8 & 3 & 8 & 3 & 4 & 6 & 8 &
1 & 8 & 1 & 6 & 4 & 9 & 2 & 9 & 2 & 4 & 5 & 7 & 3 & 7 & 3 & 5
\\
9 & 5 & 9 & 1 & 5 & 1 & 7 & 6 & 7 & 2 & 6 & 2 & 8 & 4 & 8 & 3 & 4 & 3 & 8 & 6 &
8 & 1 & 6 & 1 & 9 & 4 & 9 & 2 & 4 & 2 & 7 & 5 & 7 & 3 & 5 & 3 \end{array}
\right).\hspace*{-9pt}
\end{eqnarray*}

The design obtained from this starting design and the generating sequence
$[ 1\ 1\ 2\ 2\ 3\ 3 ]$, respectively, $[ 1\ 1\ 2\ 2\ 3\ 3\ 3 ]$,
has efficiency factor equal to $ 0.950$, respectively, to $ 0.954$.

\subsection{Comments}
Here we briefly discuss
the performances of the optimal designs obtained in this paper when the
true statistical model is simpler than the full interaction model.

Under the assumption that
the true model is the self and mixed model proposed by \citet{afsa:heda:2002}, \citet{drui:tins:2014} obtained optimal
approximate designs for the estimation of total effects. So, we can
compute the efficiency factors of our designs as defined
in (\ref{eq.defeff}) for several values of $k$ and for all $t$ with
$2\leq t\leq30$.
For $k=3$, our designs have efficiency factors greater than $0.67$.
For $k=4$, the optimal designs are the same under both models.
For $k=5$, our designs have efficiency factors greater than $0.98$. For
$k=6$, our designs have efficiency factors greater than $0.97$.

We cannot make the analogous comparison under the assumption that
the additive model is the true one, because in this case there are
no optimal designs for total effects available in the literature
[\citet{bail:drui:2004}, considered only circular designs].

We now compare our designs to complete-block neighbour-balanced designs
(CBNBDs) such as the column-complete latin squares widely used in practice.

Under the self and mixed model, CBNBDs give nonestimable total
effects but are optimal for the estimation of direct treatment effects
\citet{kune:stuf:2002}. The efficiency factors of our designs for
the direct treatment effects are $0.39$ for $k=t=3$; $0.33$ for
$k=t=4$; $0.25$ for $k=t=5$; $0.33$ for $k=t=6$; and $0.36$ for $k=t=7$.

Under the additive model, the efficiency factors of our designs for the
estimation of total effects relative to CBNBDs are $1.15$ for $k=t=3$;
$1.31$ for $k=t=4$; $1.24$ for $k=t=5$; $1.33$ for $k=t=6$; and $1.38$
for $k=t=7$. For the estimation of direct effects, CBNBDs are optimal
[\citet{kune:1984,kush:1997}], and the efficiency factors of our designs
are $0.82$ for $k=t=3$; $0.67$ for $k=t=4$; $0.52$ for $k=t=5$; $0.59$
for $k=t=6$ and $0.61$ for $k=t=7$.

\section*{Acknowledgements}
Most of this work was carried out in
the autumn of 
2011 at the Isaac Newton Institute for Mathematical Sciences in
Cambridge, UK,
during the Design and Analysis of Experiments programme.

%



\printaddresses
\end{document}